\documentclass[11pt]{article}
\usepackage{amssymb}
\usepackage{color}

\normalbaselines

\newtheorem{Theorem}{Theorem}[section]

\newtheorem{Proposition}[Theorem]{Proposition}

\newtheorem{Lemma}[Theorem]{Lemma}

\newtheorem{Corollary}[Theorem]{Corollary}

\newtheorem{Remark}[Theorem]{Remark}

\newcommand{\RR}{{{\rm I} \kern -.15em {\rm R} }}

\newcommand{\C}{{{\rm l} \kern -.42em {\rm C} }}

\newcommand{\nat}{{{\rm I} \kern -.15em {\rm N} }}

\newcommand{\be}{\begin{equation}}
\newcommand{\ee}{\end{equation}}
\newcommand{\beq}{\begin{eqnarray}}
\newcommand{\eeq}{\end{eqnarray}}
\newcommand{\beqs}{\begin{eqnarray*}}
\newcommand{\eeqs}{\end{eqnarray*}}
\newcommand{\bt}{\begin{Theorem}}
\newcommand{\et}{\end{Theorem}}
\newcommand{\br}{\begin{Remark}}
\newcommand{\er}{\end{Remark}}
\newcommand{\bc}{\begin{Corollary}}
\newcommand{\ec}{\end{Corollary}}
\newcommand{\bl}{\begin{Lemma}}
\newcommand{\el}{\end{Lemma}}
\newcommand{\bd}{\begin{definition}}
\newcommand{\ed}{\end{definition}}

\title{Stability results for  second--order evolution
 equations\\
with 
switching time--delay}                   
\author{{\sc Serge Nicaise}
\\Institut des Sciences et Techniques de Valenciennes
\\Laboratoire de Math\'ematiques et leurs Applications 
\\Universit\'e de Valenciennes et du Hainaut Cambr\'esis
\\
59313 Valenciennes Cedex 9, France\\\\
{\sc Cristina Pignotti}
\\Dipartimento di Ingegneria e Scienze dell'Informazione e Matematica\\
 Universit\`{a} di L'Aquila\\
Via Vetoio, Loc. Coppito, 67010 L'Aquila Italy}

\date{}

\begin{document}

\textwidth=160 mm

\textheight=225mm

\parindent=8mm

\frenchspacing

\maketitle

\begin{abstract}
We consider second--order evolution equations  in an abstract setting with  intermittently delayed/not--delayed damping. We give sufficient conditions for asymptotic and exponential stability,
improving and generalising
our previous results from \cite{ADE2012}.
In particular, under suitable conditions, we can consider unbounded damping 
operators. Some concrete examples         are finally presented.
\end{abstract}

\vspace{5 mm}

\def\qed{\hbox{\hskip 6pt\vrule width6pt
height7pt
depth1pt  \hskip1pt}\bigskip}



\section{Introduction}
\label{pbform}\hspace{5mm}

\setcounter{equation}{0}

Let $H$ be a real Hilbert space and let $A:{\mathcal D}(A)\rightarrow H$ 
be a positive self--adjoint  operator with a compact inverse in $H.$ Denote by $V:={\mathcal D}(A^{\frac 1 2})$ the domain of  
$A^{\frac 1 2}.$ Moreover, 
for $i=1,2,$
let $U_i$ be 
     real Hilbert spaces
with norm and inner product denoted respectively by $\Vert \cdot\Vert_{U_i}$
and $\langle\cdot ,\cdot\rangle_{U_i}$
and let
        $B_i(t):U_i\rightarrow V^\prime     ,$          be time--dependent 
linear
operators satisfying

$$B_1(t)B_2(t)=0,\quad \forall t>0.$$

Let us consider the problem
 \begin{eqnarray}
& &u_{tt}(t) +A u (t)+B_1(t)B_1^*(t)u_t(t) +B_2(t)B_2^*(t)u_t(t-\tau) =0\quad t>0\label{1.1}\\
& &u(0)=u_0\quad \mbox{\rm and}\quad u_t(0)=u_1\quad\label{1.2}
\end{eqnarray}
where the constant  $\tau >0$ is the time delay.
We assume that the delay feedback operator $B_2$ is bounded, that is
$B_2\in \mathcal{L}(U_2, H),$ while the standard one $B_1\in \mathcal{L}(U_1, V^\prime )$ may be unbounded.

Time delay
effects
appear 
in many applications and practical problems and it is by now well--known that even an arbitrarily small delay in the feedback may destabilize a system wich is uniformly exponentially stable in absence of delay.
For some examples of this destabilizing effect of time delays we refer to \cite{Datko, DLP, NPSicon06, XYL}.

In \cite{NPSicon06} we considered the wave equation with
both dampings acting simultaneously, that is $B_1(t)=\mu_1$ and $B_2(t)=\mu_2,$ 
with $\mu_1,\mu_2\in\RR^+,$ and we proved that if $\mu_1>\mu_2$ then the system is uniformly exponentially stable. Otherwise, if $\mu_2\ge\mu_1,$ that is the delay term prevails on the not delayed one, then there are instability phenomena,
namely there are unstable solutions  for arbitrarily small (large) delays.

In \cite{ADE2012}
      we firstly considered
 second--order evolution equations with intermittent delay,
this means that the standard damping and the delayed one act in different time intervals. 
This is clearly related to the stabilization problem of second--order evolution equations damped by positive/negative  feedbacks. 
We refer for this subject to \cite{HMV}. See also \cite{pignotti} for the link
between                          wave equation with time delay 
in the damping
and 
wave equation with indefinite damping, namely damping which may change sign
in different subsets of the domain.

Assuming that an observability inequality holds for the conservative model associated with (\ref{1.1})--(\ref{1.2}) and, through the definition of a suitable energy (see (\ref{energy2A})), we obtained 
in  \cite{ADE2012}
sufficient conditions ensuring asymptotic stability. Under more restrictive assumptions an exponential stability estimate was also obtained.
Our abstract framework was then applied to some concrete examples, namely the wave equation, the elasticity system and the Petrovsky system.

We mention that 
a  similar problem has been 
considered
        in \cite{ANP2012} for 1-d models for the wave equation but with a            different approach.
Indeed,
  in \cite{ANP2012} we obtain stability results for particular values of the time delays, related to the length of the domain, by using the D'Alembert formula.

Here we improve the results of \cite{ADE2012} removing a quite restrictive assumption on the size of the ``bad'' terms, that is the terms with time delay
(cfr. assumption (3.3) of \cite{ADE2012}).
Indeed, as we expect 
from \cite{HMV} and considering the relation between delay problems and 
problems with anti--damping, 
the delay feedback operator
 $B_2$ may be also large but on small time intervals.
Moreover, under an additional assumption on the size of the time intervals
where only the delay feedback acts, we give stability
results also for $B_1$ unbounded, 
while the method in \cite{ADE2012} is applicable only for bounded 
damping operators $B_1$.
Some new examples, not covered by the analysis of \cite{ADE2012}, are illustrated.
Furthermore, we point out the improvement for examples already  considered there.  

The paper is organized as follows. In section \ref{well}    a well--posedness result of the abstract system is proved. In section \ref{st} we obtain,
for $B_1$ bounded,
asymptotic and exponential stability results for the abstract model
 under suitable conditions,
improving the results of \cite{ADE2012}.
In section \ref{st2}, under an additional condition on the length of
the delay intervals, we obtain stability results valid for $B_1$
non necessarily bounded.
 Finally, in section  \ref{WW}, we illustrate our abstract results by some concrete applications.                           

\section{Well-posedness \label{well}}

\hspace{5mm}

\setcounter{equation}{0}

In this section we will give well-posedness results for problem
(\ref{1.1})--(\ref{1.2}) 
using semigroup theory. 

We assume that for all $n\in\nat$, there exists $t_n>0$ with $t_n<t_{n+1}$
and such that
\begin{eqnarray*}
B_2(t)=0\  \forall\ t\in I_{2n}=[t_{2n},t_{2n+1}),\\
B_1(t)=0 
\ \forall \  t\in I_{2n+1}=[t_{2n+1},t_{2n+2}),
\end{eqnarray*}
with  $B_2\in C([t_{2n+1},t_{2n+2}]; {\mathcal L} (U_2,  {H}))$;
for the operators $B_1$,    we assume either
$$B_1\in C^1([t_{2n},t_{2n+1}]; {\mathcal L} (U_1,  {H}))$$
or 
$$
B_1(t)=\sqrt{b_1(t)} C_n,
$$
with $C_n\in  {\mathcal L} (U_1, V'))$ and $b_1\in W^{2,\infty}(t_{2n},t_{2n+1})$
such that
$$
b_1(t)>0,\  \forall\  t\in I_{2n}.
$$

We further assume that $\tau\leq T_{2n}$ for all $n\in\nat,$ where $T_n$ denotes the lengt of the interval $I_n,$ that is
\begin{equation}\label{Tn}
T_n=t_{n+1}-t_n,\quad n\in \nat\,.
\end{equation}

Under these assumptions, we obtain the following result

\begin{Theorem}\label{texistence}
Under the above assumptions, for any $u_0 \in V$ and $u_1\in H$, the 
system $(\ref{1.1})-(\ref{1.2})$ has a unique solution 
 $u\in C([0,\infty); V)  \cap C^1([0,\infty); H)$.
\end{Theorem}
\noindent {\bf Proof.}
The case $B_1\in C^1([t_{2n},t_{2n+1}]; {\mathcal L} (U_1,  {H}))$ was treated in Theorem 2.1 of 
\cite{ADE2012}, hence we concentrate on the case when $B_1$ is not bounded.
In view of the proof of Theorem 2.1 of 
\cite{ADE2012} and the assumption on $B_1$, we only need to prove existence in the   interval $(0, t_1)$ of the problem
\begin{eqnarray}
& &u_{tt}(t) +A u (t)+ b_1(t) C_1 C_1^*u_t(t)=0\quad t>0\label{1.1bis}\\
& &u(0)=u_0\quad \mbox{\rm and}\quad u_t(0)=u_1.\quad\label{1.2bis}
\end{eqnarray}
Here the difficulty is that
a priori the domain of the operator
$$
{\cal A}(t) (u,v)^\top= (v, -Au-b_1(t) C_1 C_1^* v)^\top 
$$
will depend on $t$ since it requires that
$Au+B_1(t)B_1^*(t) v$ has to be in $H$.
Hence we cannot directly use the theory developed by Kato,  see \cite{kato:76,kato:85}.

The solution is to set (for shorthness the index 1 is suppressed in the remainder of the proof)
$$
U=(u, b u_t)^\top,
$$
that satisfies formally
$$
U_t={\cal A}(t)U,\  U(0)=(u_0, b(0) u_1)\in V\times H,
$$
where
$$
{\cal A}(t) (u,v)^\top=(b^{-1} v, b' b^{-1} v-bAu-b C C^* v).
$$
The main idea is that the domain of this new operator is independent of $t$, since it is given by
$$
D({\cal A}(t))=\{(u,v)^\top\in V\times V: Au+ C C^* v\in H\}.
$$

Now we introduce the time--depending inner product (on $V\times H$)
$$
((u,v)^\top, (\tilde u, \tilde  v)^\top)_t
=b^2(t) (A^{1/2} u, A^{1/2} \tilde u)_H+(v, \tilde v)_H, \quad \forall (u,v)^\top,\  (\tilde u, \tilde  v)^\top\in V\times H,
$$
and let $\|U\|_t=(U, U)_t^{1/2}$ its associated norm.
It is easy to check that (see for instance Theorem 2.3 of \cite{EJDE2012})
$$
\|U\|_t\leq \|U\|_s e^{\kappa |t-s|},\  \forall\  t,s\in (0,t_1),
$$
where $\kappa =\max_{t\in [0,t_1]} \frac{|b'(t)\vert }{b(t)}$.

It is a simple exercise to check that
$$
\tilde {\cal A}(t)={\cal A}(t)-\kappa Id
$$
generates a $C_0$--semigroup of contraction on $V\times H$ (that is dissipative for inner product 
$(\cdot, \cdot)_t$).
By the assumption on $b$, we deduce that $\tilde {\cal A}=\{\tilde {\cal A}(t); t\in [0,t_1]\}$
and $Y=D({\cal A}(0))$ forms a CD--system (or constant
domain system) in the sense of Kato,  see \cite{kato:76,kato:85}. In other words, for all 
$U_0=(u_0, b(0) u_1) \in V\times H$ (resp. $D({\cal A}(0))$), there exists a unique mild (resp. strong) solution $\tilde U\in C([0,t_1]; V\times H)$  (resp. $\tilde U\in C([0,t_1]; D({\cal A}(0)))\cap C^1([0,t_1]; V\times H)$)  of
$$
\tilde U_t(t)=\tilde {\cal A}(t) \tilde  U(t), t>0, \quad U(0)=U_0.
$$

Setting 
$$
U(t)=e^{\kappa t} \tilde U(t),
$$
we deduce that it is a mild (resp. strong) solution of
$$
U_t(t)=  {\cal A}(t) \tilde  U(t), t>0, \quad U(0)=U_0.
$$
Coming back to the definition of ${\cal A}(t)$, we find that $u$ is a solution of (\ref{1.1bis})--(\ref{1.2bis}).\qed

\section{Stability result: $B_1$ bounded
\label{st}}

\hspace{5mm}

\setcounter{equation}{0}

In this section we assume $B_1$ bounded. 
To get stability we assume that there exist Hilbert spaces $W_i,\ i=1,2,$ such that $H$ is continuously embedded into $W_i,$ i.e.
\begin{equation}\label{embedding}
  \|u\|_{W_i}^2\le C_i \|u\|_H^2,\quad\forall u\in H \ \mbox{\rm with}\ \  C_i>0\ \ 
\mbox{\rm  independent of}\  u.
\end{equation}
Moreover, we assume that for all $n\in\nat$, there exist three positive constants $m_{2n}$, 
 $M_{2n}$ and $M_{2n+1}$ with $m_{2n}\leq M_{2n}$ and such that for all $u\in H$ we have

i) $m_{2n}\|u\|_{W_1}^2\le  \|B_1^*(t)u\|_{U_1}
^2\le M_{2n} \|u\|_{W_1}^2$ for $t\in I_{2n}=[t_{2n},t_{2n+1}),$  $\ \forall\ n\in\nat;$

\hspace{5mm}

ii)$\|B_2^*(t)u\|_
{U_2}
^2 \le M_{2n+1}\|u\|_{W_2}^2 $ for $t\in I_{2n+1}=[t_{2n+1},t_{2n+2}),$ $\ \forall\ n\in\nat.$

\subsection{Stability without restriction on the delay time intervals
 \label{stgeneral}}

\hspace{5mm}

In this section we assume $W_1=W_2$ and we use the notation
$$W:=W_1=W_2.$$

Moreover, we assume
\begin{equation}\label{new}
\inf_{n\in\nat}
 \frac {m_{2n}}{M_{2n+1}} >0.
\end{equation}

Note that 
assumption (3.3) in
   \cite{ADE2012}   is instead equivalent to 
$$ \frac {m_{2n}}{M_{2n+1}}>
\frac 1 2.$$

Let us introduce  the energy of the system

\begin{equation}\label{energy2A}
 \displaystyle{
E(t)=E(u;t):= \frac{1}{2}\Big(\|u(t)\|_V^2+\|u_t(t)\|_H^2 +
\frac{\xi} 2
 \int_{t-\tau}^{t}\Vert B^*_2(s+\tau)u_t(s)\Vert_
{U_2}
^2  ds\Big),
}
\end{equation}
where $\xi$ is a positive number satisfying
\begin{equation}\label{suxi}
\xi <  \inf_{n\in\nat} 
\frac {m_{2n}}{M_{2n+1}}.
\end{equation}

\begin{Proposition}\label{derivE2abstrait}
Assume $\mbox{\rm i),\ ii)}$ and $(\ref{new}).$
For any regular solution of problem $(\ref{1.1})-(\ref{1.2})$ the energy is decreasing
on the intervals $I_{2n},$ $n\in\nat,$ and 
\begin{equation}\label{stimader2abstrait}
E^{\prime}(t)\le - 
\frac { m_{2n}}                {2}      \|u_t\|^2_W.
\end{equation}
Moreover, on the intervals $I_{2n+1},$ $n\in\nat,$
\begin{equation}\label{stimaderD2abstrait}
E^{\prime}(t)\le                \frac
{M _{2n+1}} 2 (\xi +\frac 1 {\xi})
\|u_t\|_W^2.
\end{equation}
\end{Proposition}

\noindent{\bf Proof:} Differentiating $E(t)$ we get
\begin{eqnarray*}
E'(t)=(u_t,u)_V+(u_{tt}, u_t)_H+
\frac{\xi}{2}\Vert B^*_2(t+\tau)u_t(t)\Vert_{U_2}^2
-\frac{\xi}{2}\Vert B^*_2(t)u_t(t-\tau)\Vert_{U_2}^2.
\end{eqnarray*}
Then,  using the definition of $A$ and (\ref{1.1}) we obtain
\begin{eqnarray*}
E'(t)&=&\langle u_t, u_{tt}+Au\rangle_{V-V'}+
\frac{\xi}{2}\Vert B^*_2(t+\tau)u_t(t)\Vert_{U_2}^2
-\frac{\xi}{2}\Vert B^*_2(t)u_t(t-\tau)\Vert_{U_2}^2\\
&=&-\langle u_t, B_1(t)B_1^*(t)u_t(t) +B_2(t)B_2^*(t)u_t(t-\tau) \rangle_{V-V'}\\
&{ }&+
\frac{\xi}{2}\Vert B^*_2(t+\tau)u_t(t)\Vert_{U_2}^2
-\frac{\xi}{2}\Vert B^*_2(t)u_t(t-\tau)\Vert_{U_2}^2.
\end{eqnarray*}
By the definition of the dual operators, we arrive at
\begin{eqnarray*}
E'(t)=-\|B_1^*(t)u_t(t)\|^2_{U_1} -(B_2^*(t)u_t,B_2^*(t)u_t(t-\tau))_{U_2}\\+
\frac{\xi}{2}\Vert B^*_2(t+\tau)u_t(t)\Vert_{U_2}^2
-\frac{\xi}{2}\Vert B^*_2(t)u_t(t-\tau)\Vert_{U_2}^2.
\end{eqnarray*}
If $t\in I_{2n}$, then $B_2(t)=0$ and the previous identity gives
\begin{eqnarray*}
E'(t)=-\|B_1^*(t)u_t(t)\|^2_{U_1} +
\frac{\xi}{2}\Vert B^*_2(t+\tau)u_t(t)\Vert_{U_2}^2.
\end{eqnarray*}
Since $T_{2n}=\vert I_{2n}\vert \geq \tau,$ it results
that $t+\tau\in I_{2n}\cup I_{2n+1}\cup I_{2n+2}.$ Now, if $t+\tau\in I_{2n}\cup I_{2n+2},$ then $B_2(t+\tau)=0.$
Therefore, $B_2(t+\tau)\ne 0$ only if $t+\tau\in I_{2n+1}.$
In both cases, 
by our assumptions i) and ii), we get (\ref{stimader2abstrait}).

For $t\in I_{2n+1},$ as $B_1(t)=0$, the previous identity becomes
\begin{eqnarray*}
E'(t)=-(B_2^*(t)u_t,B_2^*(t)u_t(t-\tau))_{U_2}+
\frac{\xi}{2}\Vert B^*_2(t+\tau)u_t(t)\Vert_{U_2}^2
-\frac{\xi}{2}\Vert B^*_2(t)u_t(t-\tau)\Vert_{U_2}^2.
\end{eqnarray*}
By Young's inequality we get
\begin{eqnarray*}
E'(t)\le\frac{1}{2\xi}\Vert B^*_2(t)u_t(t)\Vert_{U_2}^2
+\frac{\xi}{2}\Vert B^*_2(t+\tau)u_t(t)\Vert_{U_2}^2.
\end{eqnarray*}
This proves (\ref{stimaderD2abstrait}) using assumption ii) because $t+\tau$ is either in $I_{2n+1},$
or in $I_{2n+2}$ and in that last case $B^*_2(t+\tau)=0$.
$\qed$

Consider now the conservative system associated with (\ref{1.1})--(\ref{1.2})
 \begin{eqnarray}
& &w_{tt}(t) +A w (t)=0\quad t>0\label{cons1abstrait}\\
& &w(0)=w_0\quad \mbox{\rm and}\quad w_t(0)=w_1\quad\label{cons2abstrait}
\end{eqnarray}
with $(w_0,w_1)\in V\times H.$
Denote by $E_S(\cdot)$ the standard energy for wave type equations, that is
\begin{equation}\label{energystandard}
E_S(t)=E_S(w,t):=\frac 1 2( \|w\|_V^2+\|w_t\|_H^2).
\end{equation}

For our stability result we need that a suitable   observability inequality holds.
Namely we assume that there exists a time $\overline T>0$ such that for every time $T>\overline T$ there is     a constant $c,$ depending on $T$ but independent of the initial data, such that
\begin{equation}\label{OCabstrait}
E_S(0)\le c\int_0^{T}\|w_t(s)\|^2_W  ds,
\end{equation}
for every weak solution of problem $(\ref{cons1abstrait})-(\ref{cons2abstrait})$
with initial data $(w_0,w_1)\in V\times H.$

\begin{Proposition}\label{intOK2abstrait}
Assume $\mbox{\rm i), ii)}$ and $(\ref{new}).$
Moreover,
we assume that the observability inequality $(\ref{OCabstrait})$ holds for every  time $T>\overline T$
and that, denoting  $T^*:=\inf_n \{T_{2n}\},$ it is
\begin{equation}\label{quartaA}
T^*>\overline T,\quad T^*\ge \tau\,.
\end{equation}
Then, for any solution of system $(\ref{1.1})-(\ref{1.2})$ we have
\begin{equation}\label{stimabuona2A}
E(t_{2n+1})\le c_n E(t_{2n}),\quad \forall\ n\in\nat,
\end{equation}
where
\begin{equation}\label{n}
\displaystyle{
c_n= 
\frac{4c(1+4C^2T_{2n}^2M_{2n}^2) }
{m_{2n}+ 4c(1+4C^2T_{2n}^2M_{2n}^2)}
},
\end{equation}
  $c$ being the observability constant in $(\ref{OCabstrait})$ corresponding to the time $T^*$ and $C$ the constant in the norm embedding $(\ref{embedding})$  between $W$ and $H$.
\end{Proposition}

\noindent {\bf Proof.} It is sufficient to prove the estimate (\ref{stimabuona2A}) in the interval 
$I_0=[0, t_1).$ We can proceed analogously in the other intervals $I_{2n},\ n\in\nat.$

We can decompose 
$$u=w+\tilde w$$
where $w$ is a solution of system
$(\ref{cons1abstrait})-(\ref{cons2abstrait})$ with $w_0=u_0,$ $w_1=u_1;$
while $\tilde w$ solves 

\begin{eqnarray}
& &\tilde  w_{tt}(t) +A \tilde  w (t)=-B_1(t) B_1^*(t) u_t(t)\quad t>0\label{N.1babstrait}\\
& &\tilde  w(0)=0\quad \mbox{\rm and}\quad \tilde   w_t(0)=0\quad\label{N.3babstrait}
\end{eqnarray}

First we have
$$
E(0)=E_S(w,0)+ 
\frac 1 2     \int_{-\tau}^0
\| B^*_2(s+\tau) u_t(s)\|_{U_2}^2  ds=E_S(w,0),
$$
because for $s\in (-\tau,0)$, $s+\tau<\tau<t_1$.
Now using the observability inequality (\ref{OCabstrait})  we can estimate

\begin{equation}\label{E1bisabstrait}
E(0)=E_S(w,0)\le c
\int_0^{T^*}\|w_t(s)\|^2_W  ds.
\end{equation}
Using the splitting $w=u-\tilde w$ 
and the fact that $T_0=t_1\ge T^*$, we deduce that
\begin{equation}\label{E2bisabstrait}
E(0)\le 2c\int_0^{T_0}(\|{\tilde w}_t(s)\|_W^2 +\|u_t(s)\|_W^2) ds.
\end{equation}

Now, observe that from equation (\ref{N.1babstrait}) we deduce
$$
\frac{d}{dt} \frac 1 2 (\|\tilde w_t(t)\|_H^2+\|\tilde w(t)\|^2_V) =
(\tilde w_t, \tilde w_{tt}+A\tilde w)_H
= -(\tilde w_t, B_1(t) B_1^*(t) u_t(t))_H\,. 
$$
Integrating this identity  in $[0,t]$ with $0<t<T_0$, recalling (\ref{N.3babstrait}), 
and using the assumption i), we get
\begin{equation}\label{E3bisabstrait}
\begin{array}{l}
\displaystyle{
\frac 1 2 (\|\tilde w_t(t)\|_H^2+\|\tilde w(t)\|^2_V) =-\int_0^{t}(B_1^*(s) \tilde w_t(s),  B_1^*(s) u_t(s))_H  ds}
\\\medskip
\hspace{4 cm}
\displaystyle{
\le M_0\int_0^{T_0}\|\tilde w_t(s)\|_W \|u_t(s)\|_W ds.
}\end{array}
\end{equation}

Integrating (\ref{E3bisabstrait}) on $[0,t_1],$ we deduce
$$\begin{array}{l}
\displaystyle{
\int_0^{T_0} \|\tilde w_t(t)\|_W^2 dt\le C\int_0^{T_0} \|\tilde w_t(t)\|_H^2 dt
 \le 2CT_0M_0\int_0^{T_0}\|\tilde w_t(s)\|_W \|u_t(s)\|_W ds}\\
\displaystyle{
\hspace{2cm} \le CT_0M_0\int_0^{T_0}(\epsilon \|\tilde w_t(t)\|_W ^2+\frac{1}{\epsilon}\|u_t(t)\|_W^2) dt
},
\end{array}
$$
for all $\epsilon>0$
and therefore choosing $\epsilon$ such that $CT_0M_0\epsilon=\frac12$, we arrive at
\begin{equation}\label{E4bisabstrait}
\int_0^{T_0} \|\tilde w_t(t)\|_W^2 dt \le  4C^2 T_0^2M_0^2\int_0^{T_0}\|u_t(t)\|_W^2 dt\,.
\end{equation}
From (\ref{E2bisabstrait}) and  (\ref{E4bisabstrait}) 
we obtain
\begin{equation}\label{E5bisabstrait}
\begin{array}{l}
\displaystyle{
E(0)\le 2c (1+4C^2T_0^2M_0^2)
\int_0^{T_0}\|u_t(t)\|_W^2 dt}\\ 
\displaystyle{\quad
\le 
\frac {4c(1+4C^2T_0^2M_0^2)}{m_{0}}
\frac {m_0} 2 \int_0^{T_0}\|u_t(t)\|_W^2 dt.
}
\end{array}
\end{equation}
From (\ref{stimader2abstrait}) and 
(\ref{E5bisabstrait}) we deduce
$$E(t_1)\le E(0)\le 
\frac {4c(1+4C^2T_0^2M_0^2)}{m_{0}}
(E(0)-E(t_1)),$$
where we used also the fact that $E(\cdot)$ is decreasing on the time interval
$[0,t_1].$
This clearly implies 
$$E(t_1)\le c_0E(0),$$
with 
$$c_0=
\frac {4c(1+4C^2T_0^2M_0^2)}{m_{0}+4c(1+4C^2T_0^2M_0^2)}\,.\qed$$

\begin{Theorem}\label{stab2abstrait} 
Under the assumptions of Proposition $\ref{intOK2abstrait},$
if 
\begin{equation}\label{star}
\sum_{n=0}^\infty M_{2n+1} T_{2n+1}< +\infty
\quad\mbox{\rm and}\quad
       \sum_{n=0}^\infty
\frac{m_{2n}}{1+4C^2 T_{2n}^2M_{2n}^2}=+\infty\, ,    
\end{equation}
then
system $(\ref{W.1})-(\ref{W.3})$ is asymptotically stable, that is any solution  $u$ of $(\ref{W.1})-(\ref{W.3})$ satisfies 
$E(u,t)\rightarrow 0$ for $t\rightarrow +\infty\,.$ 
\end{Theorem}

\noindent {\bf Proof.} 
Note that (\ref{stimaderD2abstrait}) implies
$$E^{\prime}(t)\le M_{2n+1}(\xi +\frac 1 {\xi})CE(t),\quad t\in I_{2n+1}=[t_{2n+1},t_{2n+2}),\ n\in\nat.$$
Then we have
\begin{equation}\label{M12A}
E(t_{2n+2})\le e^{C(\xi +\frac 1 {\xi })M_{2n+1}T_{2n+1}}E(t_{2n+1}),
\quad \forall \ n\in\nat.
\end{equation} 
Combining Proposition \ref{intOK2abstrait} and (\ref{M12A}) we obtain
$$E(t_{2n+2})\le e^{C(\xi +\frac 1 {\xi })M_{2n+1} T_{2n+1}} c_n E(t_{2n}),\quad n\in\nat\,,$$
and therefore
\begin{equation}\label{M22A}
E(t_{2n+2})\le \Big (
\displaystyle{
\Pi_{p=0}^n}e^{C(\xi +\frac 1 {\xi })M_{2p+1}T_{2p+1}}c_p
\Big )E(0)\,.
\end{equation}
Then, by (\ref{M22A}), asymptotic stability occurs if
\begin{equation}\label{M32A}
  \sum_{p=0}^\infty [C(\xi +\frac 1 {\xi })M_{2p+1}T_{2p+1}+\ln c_p ]=-\infty\,.
\end{equation}
In particular  (\ref{M32A}) holds true if (\ref{star}) is valid.
Indeed, from (\ref{n}), 
$$c_p=\frac{1}
{\frac{m_{2p}}{4c(1+4C^2T_{2p}^2M_{2p}^2)}+1},$$
and then
\begin{equation}\label{zero}
\ln c_p=-\ln\Big (
1+\frac{m_{2p}}{4c(1+4C^2T_{2p}^2M_{2p}^2)}
\Big ).
\end{equation}

So, if $\frac{m_{2p}}{1+4C^2T_{2p}^2M_{2p}^2}$ tends to $0$ as $p\rightarrow\infty,$ then
$$-\ln c_p \sim \frac{m_{2p}}{4c(1+4C^2T_{2p}^2M_{2p}^2)}.$$
Consequently if (\ref{star}) holds then
$$
\sum_{p=0}^\infty \ln c_p=-\infty.
$$
Otherwise, if $\frac{m_{2p}}{1+4C^2T_{2p}^2M_{2p}^2}$ does not tend to $0,$
then, by (\ref{zero}),
$\sum_{p=0}^\infty \ln c_p=-\infty.$
Therefore, conditions (\ref{star}) imply 
(\ref{M32A}).
\qed

We now show that under additional assumptions on the coefficients $T_n, m_n, M_n$ an exponential stability  result holds.

\begin{Theorem}\label{exp}
Assume $\mbox{\rm i), ii)}$ and $(\ref{new}).$
Assume also
 that the observability inequality $(\ref{OCabstrait})$ holds for  every time
$T>\overline T$
and that
\begin{equation}\label{quartaASpecial}
T_{2n}=T^*\quad\forall\ n\in\nat ,
\end{equation}
with $T^*$ satisfying $(\ref{quartaA}),$
and
\begin{equation}\label{quartaASpecialbis}
T_{2n+1}=\tilde T\quad\forall\ n\in\nat .
\end{equation}

Moreover,  assume that
\begin{equation}\label{ASS1A}
\sup_{n\in\nat} \ e^{(\xi +\frac 1 {\xi})CM_{2n+1}\tilde T}c_n=d<1,
\end{equation}
where $c_n$ is as in $(\ref{n})$. Then, there exist two positive constants $\gamma,\mu$ 
such that
\begin{equation}\label{expestimateA}
E(t)\le \gamma e^{-\mu t} E(0),\ \ t>0,
\end{equation}
for any solution of problem  $(\ref{1.1})-(\ref{1.2}).$
\end{Theorem}

\noindent {\bf Proof.} From (\ref{ASS1A}) and (\ref{M22A})
we obtain 
$$E(T^*+\tilde T)\le d E(0),$$
and also
$$E(n(T^*+\tilde T))\le d^n E(0),\quad \forall n\in\nat.$$
Then, the energy satisfies an exponential estimate like  
(\ref{expestimateA}) (see Lemma 1 of \cite{Gugat}).\qed

\begin{Remark}{\rm
In the assumptions of Theorem \ref{exp},
from (\ref{M22A}) we can see that exponential stability also holds if instead of
(\ref{ASS1A}) we assume
$$\exists
 n\in\nat\quad \mbox{\rm such}\ \mbox{\rm that}\quad 
\Pi_{p=k(n+1)}^{k(n+1)+n} e^{(\xi +\frac 1 {\xi} )CM_{2p+1}\tilde T} c_p\le d<1,\ \ \forall\ k=0,1,2,\dots
$$
}\end{Remark}

\begin{Remark}
{\rm
Our abstract results can be applied to the examples of \cite{ADE2012},
that is damped or locally damped wave equations, elasticity system, Petrovsky 
system.
Therefore, we can improve the stability results for these models. 
}
\end{Remark}

\subsection{Stability under the restriction $T_{2n+1}\le\tau$\label{st1}}

\hspace{5mm}

We assume now that the length of the delay intervals is lower that the time 
delay, that is
\begin{equation}\label{rest}
T_{2n+1}\le\tau,\quad \forall n\in\nat\,.
\end{equation}

We look at the standard energy $E_S(\cdot).$
We can give the following estimates
on the time intervals $I_{2n}, I_{2n+1},$ $n\in\nat.$

\begin{Proposition}\label{derivE2abstraitnew}
Assume $\mbox{\rm i),\ ii)}$ and $(\ref{rest}).$
For any regular solution of problem $(\ref{1.1})-(\ref{1.2})$ the energy is decreasing
on the intervals $I_{2n},$ $n\in\nat,$ and 
\begin{equation}\label{stimader2abstraitnew}
E_S^{\prime}(t)\le - 
 m_{2n}           \|u_t(t)\|^2_{W_1}.
\end{equation}
Moreover, on the intervals $I_{2n+1},$ $n\in\nat,$
\begin{equation}\label{stimaderD2abstraitnew}
E_S^{\prime}(t)\le                
\frac {M _{2n+1}} 2 \Vert u_t(t)\Vert_{W_2}^2
+\frac  {M _{2n+1}} 2 
\|u_t(t-\tau )\|_{W_2}^2.
\end{equation}
\end{Proposition}

\noindent{\bf Proof:} Differentiating $E_S(t)$ we get
\begin{eqnarray*}
E_S'(t)=(u_t,u)_V+(u_{tt}, u_t)_H
\end{eqnarray*}
Hence using the definition of $A$ and (\ref{1.1}) we get successively
\begin{eqnarray*}
E_S'(t)&=&\langle u_t, u_{tt}+Au\rangle_{V-V'}\\
&=&-\langle u_t, B_1(t)B_1^*(t)u_t(t) +B_2(t)B_2^*(t)u_t(t-\tau) \rangle_{V-V'}.
\end{eqnarray*}
By the definition of the dual operators, we arrive at
\begin{eqnarray*}
E_S'(t)=-\|B_1^*(t)u_t(t)\|^2_{U_1} -(B_2^*(t)u_t,B_2^*(t)u_t(t-\tau))_{U_2}.
\end{eqnarray*}
If $t\in I_{2n}$, then $B_2(t)=0$ and the previous identity becomes
\begin{eqnarray*}
E_S'(t)=-\|B_1^*(t)u_t(t)\|^2_{U_1}.
\end{eqnarray*}
This gives, from i),  (\ref{stimader2abstraitnew}).

For $t\in I_{2n+1},$ as $B_1(t)=0$, the previous identity gives
\begin{eqnarray*}
E_S'(t)=-(B_2^*(t)u_t,B_2^*(t)u_t(t-\tau))_{U_2}
\le\frac{1}{2}\Vert B^*_2(t)u_t(t)\Vert_{U_2}^2
+\frac{1}{2}\Vert B^*_2(t)u_t(t-\tau)\Vert_{U_2}^2.
\end{eqnarray*}
This proves (\ref{stimaderD2abstraitnew}) using assumption ii).
$\qed$

\begin{Proposition}\label{intOK2abstraitnew}
Assume $\mbox{\rm i), ii)}$ and $(\ref{rest}).$
Moreover,
we assume that the observability inequality $(\ref{OCabstrait})$ holds for every  time $T>\overline T$
and that, denoting  $T^*:=\inf_n \{T_{2n}\},$ 
$(\ref{quartaA})$ is satisfied.
Then, for any solution of system $(\ref{1.1})-(\ref{1.2})$ we have
\begin{equation}\label{stimabuona2Anew}
E_S(t_{2n+1})\le \hat c_n E_S(t_{2n}),\quad \forall\ n\in\nat,
\end{equation}
where
\begin{equation}\label{nnew}
\displaystyle{
\hat c_n= 
\frac{2c(1+4C_1^2T_{2n}^2M_{2n}^2) }
{m_{2n}+ 2c(1+4C_1^2T_{2n}^2M_{2n}^2)}
},
\end{equation}
  $c$ being the observability constant in $(\ref{OCabstrait})$ corresponding to the time $T^*$ and $C_1$ the constant in the norm embedding $(\ref{embedding})$  between $W_1$ and $H$.
\end{Proposition}

\noindent {\bf Proof.} It is sufficient to prove the estimate 
(\ref{stimabuona2Anew}) in the interval 
$I_0=[0, t_1).$ We can proceed analogously in the other intervals $I_{2n},\ n\in\nat.$

We can decompose 
$$u=w+\tilde w$$
where $w$ is a solution of system
$(\ref{cons1abstrait})-(\ref{cons2abstrait})$ with $w_0=u_0,$ $w_1=u_1;$
while $\tilde w$ solves 

\begin{eqnarray}
& &\tilde  w_{tt}(t) +A \tilde  w (t)=-B_1(t) B_1^*(t) u_t(t)\quad t>0\label{N.1babstraitnew}\\
& &\tilde  w(0)=0\quad \mbox{\rm and}\quad \tilde   w_t(0)=0\quad\label{N.3babstraitnew}
\end{eqnarray}

First we have
$$
E_S(u,0)=E_S(w,0).$$

Then, using the observability inequality (\ref{OCabstrait}),  we can estimate

\begin{equation}\label{E1bisabstraitnew}
E_S(0)=E_S(w,0)\le c
\int_0^{T^*}\|w_t(s)\|^2_{W_1}  ds.
\end{equation}
Using the splitting $w=u-\tilde w$ 
and the fact that $T_0=t_1\ge T^*$, we deduce that
\begin{equation}\label{E2bisabstraitnew}
E_S(0)\le 2c\int_0^{T_0}(\|{\tilde w}_t(s)\|_{W_1}^2 +\|u_t(s)\|_{W_1}^2) ds.
\end{equation}

Now, observe that from equation (\ref{N.1babstrait}) we deduce
$$
\frac{d}{dt} \frac 1 2 (\|\tilde w_t(t)\|_H^2+\|\tilde w(t)\|^2_V) =
(\tilde w_t, \tilde w_{tt}+A\tilde w)_H
= -(\tilde w_t, B_1(t) B_1^*(t) u_t(t))_H\,. 
$$
Integrating this identity  in $[0,t]$ with $0<t<T_0$, recalling (\ref{N.3babstrait}), 
and using the assumption i), we get
\begin{equation}\label{E3bisabstraitnew}
\begin{array}{l}
\displaystyle{
\frac 1 2 (\|\tilde w_t(t)\|_H^2+\|\tilde w(t)\|^2_V) =-\int_0^{t}(B_1^*(s) \tilde w_t(s),  B_1^*(s) u_t(s))_H  ds}
\\\medskip
\hspace{4 cm}
\displaystyle{
\le M_0\int_0^{T_0}\|\tilde w_t(s)\|_{W_1} \|u_t(s)\|_{W_1} ds.
}\end{array}
\end{equation}

Integrating (\ref{E3bisabstraitnew}) on $[0,t_1],$ we deduce
$$\begin{array}{l}
\displaystyle{
\int_0^{T_0} \|\tilde w_t(t)\|_{W_1}^2 dt\le C_1\int_0^{T_0} \|\tilde w_t(t)\|_H^2 dt
 \le 2C_1T_0M_0\int_0^{T_0}\|\tilde w_t(s)\|_{W_1} \|u_t(s)\|_{W_1} ds}\\
\displaystyle{
\hspace{2cm} \le C_1T_0M_0\int_0^{T_0}(\epsilon \|\tilde w_t(t)\|_{W_1} ^2+\frac{1}{\epsilon}\|u_t(t)\|_{W_1}^2) dt
},
\end{array}
$$
for all $\epsilon>0$
and therefore choosing $\epsilon$ such that $C_1T_0M_0\epsilon=\frac12$, we arrive at
\begin{equation}\label{E4bisabstraitnew}
\int_0^{T_0} \|\tilde w_t(t)\|_{W_1}^2 dt \le  4C_1^2 T_0^2M_0^2\int_0^{T_0}\|u_t(t)\|_{W_1}^2 dt\,.
\end{equation}
From (\ref{E2bisabstraitnew}) and  (\ref{E4bisabstraitnew}) 
we obtain
\begin{equation}\label{E5bisabstraitnew}
\begin{array}{l}
\displaystyle{
E_S(0)\le 2c (1+4C_1^2T_0^2M_0^2)
\int_0^{T_0}\|u_t(t)\|_{W_1}^2 dt}\\ 
\displaystyle{\quad
\le 
\frac {2c(1+4C_1^2T_0^2M_0^2)}{m_{0}}
 {m_0}  \int_0^{T_0}\|u_t(t)\|_{W_1}^2 dt.
}
\end{array}
\end{equation}
From (\ref{stimader2abstraitnew}) and 
(\ref{E5bisabstraitnew}) we deduce
$$E_S(t_1)\le E_S(0)\le 
\frac {2c(1+4C_1^2T_0^2M_0^2)}{m_{0}}
(E_S(0)-E_S(t_1)),$$
where we used also the fact that $E_S(\cdot)$ is decreasing on the time interval
$[0,t_1].$
This clearly implies 
$$E_S(t_1)\le \hat c_0E_S(0),$$
with 
$$\hat c_0=
\frac {2c(1+4C_1^2T_0^2M_0^2)}{m_{0}+2c(1+4C_1^2T_0^2M_0^2)}\,.\qed$$

\begin{Theorem}\label{stab2abstraitnew} 
Under the assumptions of Proposition $\ref{intOK2abstraitnew},$
if  
$(\ref{star})$ holds,
then
system $(\ref{1.1})-(\ref{1.2})$ is asymptotically stable, that is any solution  $u$ of $(\ref{1.1})-(\ref{1.2})$ satisfies 
$E_S(u,t)\rightarrow 0$ for $t\rightarrow +\infty\,.$ 
\end{Theorem}

\noindent {\bf Proof.} 
Note that (\ref{stimaderD2abstraitnew}) implies
$$
\begin{array}{l}
E_S^{\prime}(t)\le M_{2n+1}C_2E_S(t)+ M_{2n
+1}C_2 E_S(t-\tau)\\
\hspace{1 cm}
\le M_{2n+1}C_2E_S(t)+ M_{2n
+1}C_2 E_S(t_{2n+1}) ,\quad t\in I_{2n+1}=[t_{2n+1},t_{2n+2}),\ n\in\nat,
\end{array}
$$
where we have used (\ref{rest}) and the fact that $E_S(\cdot )$ is
not increasing in the time intervals $I_{2n}.$
Remark that the constant $C_2$ is the one from the norm embedding (\ref{embedding})
between $W_2$ and $H.$

Then we have
\begin{equation}\label{M12Anew}
E_S(t)\le e^{M_{2n+1}C_2 (t-t_{2n+1})} E_S(t_{2n+1})+\left [
 e^{M_{2n+1}C_2 (t-t_{2n+1})}-1
\right ]E_S(t_{2n}),
\end{equation} 
for $t\in I_{2n+1}=[t_{2n+1},t_{2n+2}),\ n\in\nat.$
Combining Proposition \ref{intOK2abstraitnew} and (\ref{M12Anew}) we obtain
$$E_S(t_{2n+2})\le \left [e^{M_{2n+1} T_{2n+1}C_2} \hat c_n +e^{M_{2n+1} T_{2n+1}C_2}-1
\right ]
E_S(t_{2n}),\quad n\in\nat\,,$$
and therefore
\begin{equation}\label{M22Anew}
E_S(t_{2n+2})\le \Big (
\displaystyle{
\Pi_{p=0}^n}\left [e^{M_{2p+1} T_{2p+1}C_2} \hat c_p +e^{M_{2p+1} T_{2p+1}C_2}-1 \right ]
\Big )E_S(0)\,.
\end{equation}
Then, by (\ref{M22Anew}), asymptotic stability occurs if
\begin{equation}\label{M32Anew}
  \sum_{p=0}^\infty \left [C_2M_{2p+1}T_{2p+1}+\ln \left (\hat c_p+1-
e^{-M_{2p+1}C_2T_{2p+1}}\right )
\right
 ]=-\infty\,.
\end{equation}
In particular  (\ref{M32Anew}) holds true if (\ref{star}) is valid.
Indeed, if (\ref{star}) holds, then
$$\lim_{n\rightarrow\infty} M_{2n+1}T_{2n+1}=0,$$
and therefore, being
from (\ref{nnew}), 
$$\hat c_p=\frac{1}
{\frac{m_{2p}}{2c(1+4C_1^2T_{2p}^2M_{2p}^2)}+1},$$
it results
\begin{equation}\label{zeronew}
\ln \left (\hat c_p 
+1-
e^{-M_{2p+1}C_2T_{2p+1}}\right )\sim \ln\hat c_p
=-\ln\Big (
1+\frac{m_{2p}}{2c(1+4C_1^2T_{2p}^2M_{2p}^2)}
\Big ).
\end{equation}

So, if $\frac{m_{2p}}{1+4C_1^2T_{2p}^2M_{2p}^2}$ tends to $0$ as $p\rightarrow\infty,$ then
$$-\ln \hat c_p \sim \frac{m_{2p}}{2c(1+4C_1^2T_{2p}^2M_{2p}^2)}.$$

Otherwise, if $\frac{m_{2p}}{1+4C_1^2T_{2p}^2M_{2p}^2}$ does not tend to $0,$
then, by (\ref{zeronew}),
$\sum_{p=0}^\infty \ln \hat c_p=-\infty.$
Therefore, conditions (\ref{star}) imply
$$
\sum_{p=0}^\infty \ln \hat c_p=-\infty
$$ 
and then
(\ref{M32Anew}).
\qed

Also in this case, under additional assumptions on the coefficients $T_n, m_n, M_n,$ an exponential stability  result holds.

\begin{Theorem}\label{expnew}
Assume $\mbox{\rm i), ii)}$ and $(\ref{rest}).$
Assume also
 that the observability inequality $(\ref{OCabstrait})$ holds for  every time
$T>\overline T$
and that
\begin{equation}\label{quartaASpecialnew}
T_{2n}=T^*\quad\forall\ n\in\nat ,
\end{equation}
with $T^*$ satisfying $(\ref{quartaA}),$
and
\begin{equation}\label{quartaASpecialnewbis}
T_{2n+1}=\tilde T\quad\forall\ n\in\nat ,
\end{equation}
with $\tilde T\le\tau.$

Moreover,  assume that
\begin{equation}\label{ASS1Anew}
\sup_{n\in\nat} \left [ 
e^{C_2M_{2n+1}\tilde T}(\hat c_n+1)-1
\right ]
 =\hat d<1,
\end{equation}
where $\hat c_n$ is as in $(\ref{nnew})$. Then, there exist two positive constants $\hat\gamma,\hat \mu$ 
such that
\begin{equation}\label{expestimateAnew}
E_S(t)\le \hat\gamma e^{-\hat \mu t} E_S(0),\ \ t>0,
\end{equation}
for any solution of problem  $(\ref{1.1})-(\ref{1.2}).$
\end{Theorem}

\begin{Remark}{\rm
In the assumptions of Theorem \ref{expnew},
from (\ref{M22Anew}) we can see that exponential stability also holds if instead of
(\ref{ASS1Anew}) we assume
$$\exists
 n\in\nat\quad \mbox{\rm such}\ \mbox{\rm that}\quad 
\Pi_{p=k(n+1)}^{k(n+1)+n}\left [ 
e^{C_2M_{2p+1}\tilde T}(\hat c_p+1)-1
\right ] \le \hat d<1,\ \ \forall\ k=0,1,2,\dots
$$
}\end{Remark}

\begin{Remark}
{\rm
Our abstract results can be applied to the examples of \cite{ADE2012},
that is damped or locally damped wave equations, elasticity system, Petrovsky 
system when $T_{2n+1}\le\tau,$ $\forall\ n\in\nat.$
}
\end{Remark}

\section{Stability result: $B_1$ unbounded
\label{st2}}

\hspace{5mm}

\setcounter{equation}{0}

In this section $B_1$ may be unbounded. 
We assume that there exists a Hilbert space $W$ such that $H$ is continuously embedded into $W,$ i.e.,
\begin{equation}\label{embeddingbis}
  \|u\|_{W}^2\le C \|u\|_H^2,\quad\forall u\in H \ \mbox{\rm with}\ \  C>0\ \ 
\mbox{\rm  independent of}\  u.
\end{equation}

Moreover, we assume that $V$ is embedded into $U_1$ and that for all $n\in\nat$, there exist three positive constants $m_{2n}$, 
 $M_{2n}$ and $M_{2n+1}$ with $m_{2n}\leq M_{2n}$ such that 

i) $m_{2n}\|u\|_{U_1}^2\le  \|B_1^*(t)u\|_{U_1}
^2\le M_{2n} \|u\|_{U_1}^2$ for $t\in I_{2n}=[t_{2n},t_{2n+1}),$  
$\ \forall u\in V,\ $
$\ \forall\ n\in\nat;$

\hspace{5mm}

ii) $\|B_2^*(t)u\|_
{U_2}
^2 \le M_{2n+1}\|u\|_W^2 $ for $t\in I_{2n+1}=[t_{2n+1},t_{2n+2}),$ 
$\ \forall u\in H,\ $ 
$\ \forall\ n\in\nat.$

In order to deal with unbounded feedback we will work with the standard 
energy $E_S(\cdot ).$ 
Then, as before, we assume (\ref{rest}).

As before we can give the following estimates
on the time intervals $I_{2n}, I_{2n+1},$ $n\in\nat.$

\begin{Proposition}\label{derivE2abstraitnewU}
Assume $\mbox{\rm i),\ ii)}$ and $(\ref{rest}).$
For any regular solution of problem $(\ref{1.1})-(\ref{1.2})$ the energy is decreasing
on the intervals $I_{2n},$ $n\in\nat,$ and 
\begin{equation}\label{stimader2abstraitnewU}
E_S^{\prime}(t)\le - 
            \|B_1^*(t)u_t(t)\|^2_{U_1}.
\end{equation}
Moreover, on the intervals $I_{2n+1},$ $n\in\nat,$
the estimate $(\ref{stimaderD2abstraitnew})$ holds (with $W$ instead of $W_2.$)
\end{Proposition}

Consider now the damped system
\begin{eqnarray}
& &w_{tt}(t) +A w (t)+B_1(t) B_1^*(t) w_t=0,\quad t\in (t_{2n}, t_{2n+1}),\ n\in\nat, \label{cons1abstraitU}\\
& &w(t_{2n})=w_0^n\quad \mbox{\rm and}\quad w_t(t_{2n})=w_1^n\quad\label{cons2abstraitU}
\end{eqnarray}
with $(w_0^n,w_1^n)\in V\times H.$
For our stability result we need that the next    observability type  inequality holds.
Namely we assume that, for every $n$  there exists a time $\overline T_n,$ 
such that
\begin{equation}\label{T2n}
T_{2n}>\overline T_n,
\end{equation}
and
for every $n$ and every time $T,$ with $T_{2n}\ge T> \overline T_n,$
 there is     a constant $d_n,$ depending on $T$ but independent of $(w_0^n,w_1^n),$ such that
\begin{equation}\label{OCabstraitCris}
E_S(t_{2n}+T)\le d_n\int_{t_{2n}}^{t_{2n}+T}\Vert B_1^*(t) w_t(t)\Vert_{U_1}^2  dt,
\end{equation}
for every weak solution of problem $(\ref{1.1}), (\ref{1.2})$
with initial data $(w_0^n,w_1^n)\in V\times H.$

\begin{Proposition}\label{intOK2abstraitnewU}
Assume $\mbox{\rm i), ii)},$$(\ref{rest})$ and $T_{2n}\ge\tau,$ $\forall\ n\in\nat.$
Moreover, we assume that there is a sequence $\{\overline T_n\}_n,$ such that  
$(\ref{T2n})$ is satisfied and the
inequality $(\ref{OCabstraitCris})$ holds for every $T\in (\overline T_n, T_{2n}],$ $\forall\ n\in\nat.$ 
Then, for any solution of system $(\ref{1.1})-(\ref{1.2})$ we have
\begin{equation}\label{stimabuona2AnewU}
E_S(t_{2n+1})\le \hat d_n E_S(t_{2n}),\quad \forall\ n\in\nat,
\end{equation}
where
\begin{equation}\label{nnewU}
\displaystyle{
\hat d_n= 
\frac{d_n}{d_{n}+1}
},
\end{equation}
$d_n$ being the observability constant in 
$(\ref{OCabstraitCris})$ corresponding to the time $T_{2n}.$
 \end{Proposition}

\noindent {\bf Proof.} To prove $(\ref {stimabuona2AnewU})$ is sufficient to 
use the estimate $(\ref{stimader2abstraitnewU})$ in
$(\ref{OCabstraitCris}),$ reminding that $B_2(t)=0$ on $(t_{2n},t_{2n+1})$. $\qed$

\begin{Theorem}\label{stab2abstraitnewU} 
Under the assumptions and with the same notations of Proposition $\ref{intOK2abstraitnewU},$
if  
\begin{equation}\label{starstar}
\sum_{n=0}^\infty  M_{2n+1}T_{2n+1}<+\infty \quad
\mbox{\rm and}\quad \sum_{n=0}^\infty\ln \hat d_n =-\infty,
\end{equation}
then
system $(\ref{1.1})-(\ref{1.2})$ is asymptotically stable, that is any solution  $u$ of $(\ref{1.1})-(\ref{1.2})$ satisfies 
$E_S(u,t)\rightarrow 0$ for $t\rightarrow +\infty\,.$ 
\end{Theorem}

\begin{Remark}\label{rkserge1}
{\rm
In  fact $d_n$ depends on  $n$ because by hypothesis  $B_1$ may depend on the time variable. However,
if $B_1$ does not depend on $t,$ then by a translation of $t_{2n}$
the constant   $d_n$ becomes independent of $n$. But if $d_n=d>0$ for all $n$, then the condition
$$
\sum_{n=0}^\infty\ln \hat d_n =-\infty
$$
is automatically satisfied. On the other hand, the first condition in (\ref{starstar}) depends only on the length of the intervals $I_{2n+1}$ and on the boundedness constant of $B_2^*$ on the same intervals,
hence (\ref{starstar}) can be easily checked.
}
\end{Remark}

Also in this case, under additional assumptions on the coefficients $T_n, m_n, M_n,$ an exponential stability  result holds.

\begin{Theorem}\label{expnewU}
Assume $\mbox{\rm i), ii)}$ and $(\ref{rest}).$
Assume also that $(\ref{quartaASpecialnew})$ holds
with $T^*$ satisfying   $T^*\geq\tau$ and
that  inequality $(\ref{OCabstraitCris})$  holds, $\forall\ n\in\nat,$ for  every time
$T$ with $T^*\ge T>\overline T.$
Moreover,  assume  $T_{2n+1}=\tilde T,$ for all $ n\in\nat,$  with $\tilde T
\le \tau.$
If 
\begin{equation}\label{ASS1AnewUU}
\sup_{n\in\nat} \left [ 
e^{CM_{2n+1}\tilde T}(\hat d_n+1)-1
\right ]
 <1,
\end{equation}
where $\hat d_n$ is as in $(\ref{nnewU}),$ then, there exist two positive constants $\hat\gamma,\hat \mu$ 
such that
\begin{equation}\label{expestimateAnewU}
E_S(t)\le \hat\gamma e^{-\hat \mu t} E_S(0),\ \ t>0,
\end{equation}
for any solution of problem  $(\ref{1.1})-(\ref{1.2}).$
\end{Theorem}

\begin{Remark}\label{rkserge2}
{\rm
If  $\hat d_n\leq \hat d<1$ (see Remark \ref{rkserge1}), then  (\ref{ASS1AnewUU}) holds if
$$
(1+\hat d)\sup_{n\in\nat} \left [ 
e^{CM_{2n+1}\tilde T}-\frac 1 {1+\hat d}
\right ]
 <1,
$$
or equivalently
$$
\sup_{n\in\nat} \left [ 
e^{CM_{2n+1}\tilde T}
\right ]
 <\frac{2}{\hat d+1}.
$$
Hence (\ref{ASS1AnewUU}) is verified  if $\sup_{n\in\nat} M_{2n+1}$ is small enough.
This is a quite realistic assumption because then the influence of the delay term is small and the action of the standard dissipation  sufficiently compensates it to guarantee an exponential decay.
}
\end{Remark}

\section{Examples\label{WW}}

\hspace{5mm}

Here we apply our abstract results to some concrete models. Note that 
the following examples are not included in the setting of \cite{ADE2012}.

\setcounter{equation}{0}

\subsection{The wave equation with internal and boundary dampings}

Our first application concerns the wave equation with boundary feedback 
and internal delay term.
Let $\Omega\subset\RR^n$ be an open bounded domain   with a boundary
$\partial\Omega$ of class $C^2.$
We assume that $\partial\Omega$ is
composed of two closed sets $\partial\Omega =\Gamma_0\cup \Gamma_1,$ with
$\Gamma_0\cap\Gamma_1=\emptyset$ and $\mbox{\rm meas}\ \Gamma_1 >0.$  
  
Denoting by $m$ the standard multiplier 
$m(x)=x-x_0,\ x_0\in\RR^n,$
we assume that the  
\be\label{gamma1}
m(x)\cdot \nu (x)\le 0,\quad \mbox{\rm for}\ x\in\Gamma_1, 
\ee
and, for some $\delta >0,$
\be\label{gamma0}
m(x)\cdot \nu (x)\ge \delta,\quad \mbox{\rm for}\ x\in\Gamma_0, 
\ee
where $\nu(x)$ is the outer unit normal vector at $x\in\partial\Omega.$
Given $\omega \subseteq\Omega,$ 
let us consider the initial boundary value  problem

 \begin{eqnarray}
& &u_{tt}(x,t) -\Delta u (x,t)+b_2(t)\chi_\omega u_t(x,t-\tau)=0\quad \mbox{\rm in}\quad\Omega\times
(0,+\infty)\label{W.1}\\
& &u (x,t) =0\quad \mbox{\rm on}\quad \Gamma_1\times
(0,+\infty)\label{W.2}\\
& &\frac {\partial u}{\partial\nu} (x,t)=-b_1(t)u_t(x,t)\quad \mbox{\rm on}\quad \Gamma_0\times
(0,+\infty)\label{W.2B}\\ 
& &u(x,0)=u_0(x)\quad \mbox{\rm and}\quad u_t(x,0)=u_1(x)\quad \hbox{\rm
in}\quad\Omega\label{W.3}
\end{eqnarray}
with
initial
data $(u_0, u_1)\in H^1_{\Gamma_1}(\Omega)\times L^2(\Omega),$  
where as usual 
$$H^1_{\Gamma_1}(\Omega ):=\{\ u\in H^1(\Omega)\,:\, u=0\ \mbox{\rm on}\ \Gamma_1\,\},$$
and
$b_1,b_2$ in $L^\infty (0,+\infty).$ 

On the feedback functions $b_1(\cdot ), b_2(\cdot ),$ 
we assume 
$$b_1(t)b_2(t)=0,\quad\forall\ t>0,$$
in order to have an intermittent delay problem.
We refer to \cite{ANP2010} 
for the analysis of this problem when  $b_1, b_2$ are constant in time,
in other words the delayed damping and the standard boundary one are acting
simultaneously for every time $t>0.$

Moreover, we assume $b_1\in W^{2,\infty}(I_{2n}),$ $\forall\ n\in\nat,$ and 

$i_w$)  $0< m_{2n}\le  b_1(t)\le M_{2n}$,
$b_2(t)=0,$
for all $t\in I_{2n}=[t_{2n},t_{2n+1}),$ and $b_1\in C^1(\bar I_{2n})$, for all $\ n\in\nat;$

 $ii_w$) $\vert b_2(t)\vert \le M_{2n+1}$,
  $b_1(t)=0,$
for all $t\in I_{2n+1}=[t_{2n+1},t_{2n+2}),$ and $b_2\in C(\bar I_{2n+1})$, for all  $\ n\in\nat.$

This problem enters into our previous framework, if we take
$H=L^2(\Omega)$ and the operator  $A$ defined by
$$A:{\mathcal D}(A)\rightarrow H\,:\,  u\rightarrow -\Delta u,$$
where 
$${\mathcal D}(A)=:\{\ u\in H^1_{\Gamma_1}(\Omega)\,:\, \Delta u\in L^2(\Omega )
\ \mbox{\rm and}\ \frac {\partial u}{\partial\nu }=0 
\ \mbox{\rm on}\ \Gamma_0\,\}.$$

We then define $U_1:=L^2(\Gamma_0),$ $U_2:=L^2(\omega)$ 
and the operators $B_1(t), B_2(t )$ as
$$B_2 \in {\cal L}(U_2; H), \, B_2u  = \sqrt{b_2} \, \tilde u,
\  \forall \, u \in L^2(\omega),$$
and
$$B_1 \in {\cal L}(U_1;V^\prime    ), \, B_1u = \sqrt{b_1} \, A_{-1} Nu, \, \forall \, u \in L^2(\Gamma_0), \,
B_1^*w =\sqrt{b_1}  w_{|\Gamma_0}, \, \forall \, w \in 
V:=
{\cal D}(A^{1/2}),$$          where                                                                      $A_{-1}$ is the extension of $A$ to $H$,
namely for all $h\in H$ and $\varphi \in {\cal D}(A)$, $A_{-1}h$ is
the unique element in $({\mathcal D}(A))^\prime$ 
                                (the duality is
in the sense of $H$),
      such that (see for instance
\cite{TucsnakWeiss})
\[
\langle A_{-1}h; \varphi\rangle_{({\mathcal D}(A))^\prime,  {\mathcal D}(A)}                  =\int_\Omega h
A\varphi\, dx.\]

Here and below $N \in {\cal L}(L^2(\Gamma_0);L^2(\Omega))$ is defined as follows:  for all
$v \in L^2(\Gamma_0), \, Nv$ is the unique solution (transposition
solution) of
$$\Delta Nv =0, \,
Nv_{|\Gamma_1} = 0, \,
\frac{\partial Nv}{\partial \nu}_{|\Gamma_0} = v.$$

With these definitions, we can show that problem (\ref{W.1})--(\ref{W.3}) enters in the abstract framework
(\ref{1.1})--(\ref{1.2}) and that the assumptions i) and ii) of section \ref{st2} hold with $W=L^2(\omega)$.

As $B_1$ is not bounded, we need to consider the non delayed  system

 \begin{eqnarray}
& &w_{tt}(x,t) -\Delta w (x,t)=0\quad \mbox{\rm in}\quad\Omega\times
(0,+\infty)\label{C.1B}\\
& &w (x,t) =0\quad \mbox{\rm on}\quad\Gamma_1\times
(0,+\infty)\label{C.2B}\\
& &\frac{\partial w}{\partial \nu}(x,t)=- f(t) w_t(x,t), \quad x\in  \Gamma_0,\ t> 0\label{a3B}\\
& &w(x,0)=w_0(x)\quad \mbox{\rm and}\quad w_t(x,0)=w_1(x)\quad \hbox{\rm
in}\quad\Omega\label{C.3B}
\end{eqnarray}
with $(w_0,w_1)\in H_{\Gamma_1}^1(\Omega)\times L^2(\Omega)$
and $f\in L^\infty (0,+\infty),$  $f(t)\ge 0$ a.e. $t>0.$

\begin{Proposition}\label{KomZua}
There exists a time $\overline T>0$ such that for every $T>\overline T,$
there are constants $\alpha_i, i=1,2,3,$ for which
\begin{equation}\label{quasiob}
(T-\overline T) E_S(T)\le\alpha_1
\int_0^T\int_{\Gamma_0}\left (\frac{\partial w}{\partial \nu}\right )^2(x,t) d\Gamma dt +\alpha_2
\int_0^T\int_{\Gamma_0}f(t)w_t^2(x,t) d\Gamma dt +
\alpha_3
\int_0^T\int_{\Gamma_0}w_t^2(x,t) d\Gamma dt,
\end{equation}
for any weak solution of $(\ref{C.1B})-(\ref{C.3B}).$
The constants $\alpha_i, i=1,2,3,$ are independent of the initial data 
and of the function $f(\cdot ),$
but they 
depend on $T$ and on $\Omega.$
\end{Proposition}

\noindent {\bf Proof.} 
The estimate (\ref{quasiob}) can be easily obtained from a standard multiplier argument (cfr. (3.11)--(3.16) of \cite{KomornikZuazua}).\qed

Then, from Proposition \ref{KomZua}, we deduce that
there exists a time $\overline T>0$  such that, assuming $T_{2n}>\overline T,$
$\forall\ n\in\nat,$ then for every $n$ and every time $T$ with $T_{2n}\ge T>
\overline T,$ there is a constant $d_n$ for which  
$(\ref{OCabstraitCris})$ holds for any weak solution of $(\ref{C.1B})-(\ref{C.3B}).$

From $(\ref{quasiob})$ with $f(t)=b_1(t-t_{2n})$ and the boundary condition (\ref{W.2B}) we deduce
the explicit dependence of $d_n$ from the feedback function $b_1,$ 
that is

\begin{equation}\label{explicit}
E_S(t_{2n+1})\le d_n
\int_{t_{2n}}^{t_{2n+1}}\int_{\Gamma_0}b_1(t) u_t^2(x,t) dt d\Gamma,\ \forall\ n\in\nat,
\end{equation}
with 
\begin{equation}\label{dnexp}
d_n:=
\frac{\alpha_1M_{2n}m_{2n}+\alpha_2m_{2n}+\alpha_3}
{m_{2n}(T_{2n}-\overline T)},\quad \forall\ n\in\nat.
\end{equation}

Therefore, we can restate Theorem \ref{stab2abstraitnewU} 
under more explicit conditions.

\begin{Theorem}\label{Stabexplicit} 
Under the assumptions $i_w), ii_w)$ and $(\ref{rest}),$
if  $T^*:=\inf_n T_{2n}$ satisfies
$(\ref{quartaA})$ and 
\begin{equation}\label{starstarexplicit}
\sum_{n=0}^\infty  M_{2n+1}T_{2n+1}<+\infty \quad
\mbox{\rm and}\quad \sum_{n=0}^\infty
\frac {m_{2n}(T_{2n}-\overline T)}
{\alpha_1M_{2n}m_{2n}+\alpha_2m_{2n}+\alpha_3}
=+\infty,
\end{equation}
then
system $(\ref{W.1})-(\ref{W.3})$ is asymptotically stable, that is any solution  $u$ of $(\ref{W.1})-(\ref{W.3})$ satisfies 
$E_S(u,t)\rightarrow 0$ for $t\rightarrow +\infty\,.$ 
\end{Theorem}

\noindent{\bf Proof.}
We have only to show that the second condition of (\ref{starstarexplicit}) implies
the second condition of (\ref{starstar}), namely in this case
$$\sum_{n=0}^\infty \ln \frac {\alpha_1M_{2n}m_{2n}+\alpha_2m_{2n}+\alpha_3}
{\alpha_1M_{2n}m_{2n}+\alpha_2m_{2n}+\alpha_3+ m_{2n}(T_{2n}-\overline T)}=-\infty.$$
Now, observe that 
$$\ln\frac {\alpha_1M_{2n}m_{2n}+\alpha_2m_{2n}+\alpha_3}
{\alpha_1M_{2n}m_{2n}+\alpha_2m_{2n}+\alpha_3+ m_{2n}(T_{2n}-\overline T)}
=-\ln \left (
1+
\frac {m_{2n}(T_{2n}-\overline T)}
{\alpha_1M_{2n}m_{2n}+\alpha_2m_{2n}+\alpha_3}
\right ),
$$
then we can conclude arguing as in the proof of Theorem \ref{stab2abstraitnew}.
\qed

\begin{Remark}\label{rkserge3}
{\rm
If $T_{2n}=T^*$ and $M_{2n}=M^*$ as well as $m_{2n}=m^*$, then
 (\ref{starstarexplicit}) holds if the easily checked condition
$$
\sum_{n=0}^\infty  M_{2n+1}T_{2n+1}<+\infty
$$
holds (since 
 the second condition  of (\ref{starstarexplicit}) automatically  holds).
}
\end{Remark}

Similarly using Theorem \ref{expnewU}, we directly can state the

\begin{Theorem}\label{expnewUexplicit}
Assume that $i_w), ii_w)$ and $(\ref{rest})$ hold,
 that $T_{2n}=T^*$, for all $n\in\nat$,
with $T^*$ satisfying $T^*\geq\tau$ and $T^*>\bar T$ with $\bar T$ from Proposition $\ref{KomZua}.$
Moreover,  assume  $T_{2n+1}=\tilde T,$ for all $n\in\nat$ with $\tilde T
\le \tau.$
If 
\begin{equation}\label{ASS1AnewU}
\sup_{n\in\nat} \left [ 
e^{CM_{2n+1}\tilde T}(\hat d_n+1)-1
\right ]
 <1,
\end{equation}
where 
$$\hat d_n=\frac {\alpha_1M_{2n}m_{2n}+\alpha_2m_{2n}+\alpha_3}
{\alpha_1M_{2n}m_{2n}+\alpha_2m_{2n}+\alpha_3+ m_{2n}(T_{2n}-\overline T)},
$$ then there exist two positive constants $\hat\gamma,\hat \mu$ 
such that
\begin{equation}\label{expestimateAnewUU}
E_S(t)\le \hat\gamma e^{-\hat \mu t} E_S(0),\ \ t>0,
\end{equation}
for any solution of problem  $(\ref{W.1})-(\ref{W.3}).$
\end{Theorem}

As in the abstract setting (see Remarks \ref{rkserge1} and \ref{rkserge2}),  explicit conditions on $b_1, b_2$ and $T_{2n}$ can be found in order to get exponential decay.

\subsection{The wave equation with internal delayed/undelayed feedbacks}

Here we consider the wave equation with local internal damping
and internal delay.
More precisely, let $\Omega\subset\RR^n$ be an open bounded domain   with a boundary
$\partial\Omega$ of class $C^2.$  
Denoting by $m,$ as before, the standard multiplier 
$m(x)=x-x_0,\ x_0\in\RR^n,$ 
let $\omega_1$ be the intersection of $\Omega$ with an open neighborhood of the
subset of $\partial\Omega$
\be\label{defgamma0}\Gamma_0=\{\, x\in\partial\Omega\, :\ m(x)\cdot \nu (x)>0\,\}.
\ee

Let us consider the initial boundary value  problem

 \begin{eqnarray}
& &u_{tt}(x,t) -\Delta u (x,t)+b_1(t)\chi_{\omega_1} u_t(x,t)+b_2(t)\chi_{\omega_2} u_t(x,t-\tau)=0\ \mbox{\rm in}\ \Omega\times
(0,+\infty)\label{Wo.1}\\
& &u (x,t) =0\quad \mbox{\rm on}\quad\partial\Omega\times
(0,+\infty)\label{Wo.2}\\
& &u(x,0)=u_0(x)\quad \mbox{\rm and}\quad u_t(x,0)=u_1(x)\quad \hbox{\rm
in}\quad\Omega\label{Wo.3}
\end{eqnarray}
with
initial
data $(u_0, u_1)\in H^1_0(\Omega)\times L^2(\Omega)$  
and
$b_1,b_2$ in $L^\infty (0,+\infty)$ such that
$$b_1(t)b_2(t)=0,\quad\forall\ t>0.$$

Moreover, we assume
$i_w)$ and 
 $ii_w).$

This problem enters into our previous framework, if we take
$H=L^2(\Omega)$ and the operator  $A$ defined by
$$A:{\mathcal D}(A)\rightarrow H\,:\,  u\rightarrow -\Delta u,$$
where ${\mathcal D}(A)=H^1_0(\Omega)\cap
H^2(\Omega).$  

The operator $A$ is a self--adjoint and positive operator with a compact inverse in $H$
and is such that $V={\mathcal D}(A^{1/2})=H^1_0(\Omega).$
We then define $U_i=L^2(\omega_i)$ and the operators $B_i, i=1,2,$ as
\be\label{defBi}
B_i:U_i\rightarrow H: \quad v\rightarrow \sqrt{b_i}\tilde v \chi_{\omega_i},
\ee
where $\tilde v\in L^2(\Omega)$ is the extension of $v$ by zero outside $\omega_i.$ 
It is easy to verify that 
$$B_i^*(\varphi )=\sqrt {b_i} \varphi_{\vert_{\omega_i}}\quad\mbox{\rm for}\ \varphi\in H,$$
and thus 
(\ref{embedding}) holds with $W_i=L^2(\omega_i)$, while
$B_iB_i^*(\varphi)=b_i\varphi \chi_{\omega_i},$ for $\varphi\in H$ and $i=1,2.$
This shows that problem (\ref{Wo.1})--(\ref{Wo.3}) enters in the abstract framework
(\ref{1.1})--(\ref{1.2}).
Moreover, $i_w)$ and $ii_w)$ easily imply $i)$ and $ii)$ of sect. \ref{st}.
Therefore we can restate Proposition \ref{derivE2abstrait}.
Now, the energy functional is

\begin{equation}\label{energy2o}
 \displaystyle{
E(t)= \frac{1}{2}\int_\Omega \{ u_t^2(x,t) +\vert \nabla u(x,t)\vert^2\} dx +
 \frac \xi 2\int_{t-\tau}^{t}\vert b_2(s+\tau)\vert  \int_{\omega} u_t^2(x,s) dx ds.
}
\end{equation}

Consider now the conservative system
\begin{eqnarray} 
& &w_{tt}(x,t) -\Delta w (x,t)=0\quad \mbox{\rm in}\quad\Omega\times
(0,+\infty)\label{C.1o}\\
& &w (x,t) =0\quad \mbox{\rm on}\quad\partial\Omega\times
(0,+\infty)\label{C.2o}\\
& &w(x,0)=w_0(x)\quad \mbox{\rm and}\quad w_t(x,0)=w_1(x)\quad \hbox{\rm
in}\quad\Omega\label{C.3o}
\end{eqnarray}
with $(w_0,w_1)\in H_0^1(\Omega)\times L^2(\Omega).$
It is well--known that an observability inequality holds (see e.g. \cite{BLR, Komornikbook, KomLoretti, Lag83, LT, Lions, liu, zuazua}): 
There exists a time $\overline T>0$ such that for every time $T>\overline T$ 
there is a
constant $c,$ depending on $T$ but independent of the initial data, such that
\begin{equation}\label{OCo}
E_S(0)\le c\int_0^{T}\int_{\omega_1}w_t^2(x,s) dx ds,
\end{equation}
for every weak solution of problem $(\ref{C.1o})-(\ref{C.3o}).$

In the case $\omega_1=\omega_2$  we can apply the results of section \ref{stgeneral}. Therefore we can restate Proposition \ref{intOK2abstrait}
and Theorems \ref{stab2abstrait} and \ref{exp}.

\begin{Proposition}\label{intOK2abstraito}
Assume $\omega_1 =\omega_2$,
 $\mbox{\rm i}_w), \mbox{\rm ii}_w)$ and $(\ref{new})$ are satisfied.
Moreover,
we assume that the observability inequality $(\ref{OCo})$ holds for every  time $T>\overline T$
and that, denoting  $T^*:=\inf_n \{T_{2n}\},$  the assumption $(\ref{quartaA})$ holds.
Then, for any solution of system $(\ref{Wo.1})-(\ref{Wo.3})$ we have
\begin{equation}\label{stimabuona2Ao}
E(t_{2n+1})\le c_n E(t_{2n}),\quad \forall\ n\in\nat,
\end{equation}
where
\begin{equation}\label{no}
\displaystyle{
c_n= 
\frac{4c(1+4T_{2n}^2M_{2n}^2) }
{m_{2n}+ 4c(1+4T_{2n}^2M_{2n}^2)}
},
\end{equation}
  $c$ being the observability constant in $(\ref{OCo})$ corresponding to the time $T^*$ and $C$ the constant in the norm embedding $(\ref{embedding})$  between $W$ and $H$.
\end{Proposition}

\begin{Theorem}\label{stab2abstraito} 
Under the assumptions of Proposition $\ref{intOK2abstraito},$
if 
\begin{equation}\label{staro}
\sum_{n=0}^\infty M_{2n+1} T_{2n+1}< +\infty
\quad\mbox{\rm and}\quad
       \sum_{n=0}^\infty
\frac{m_{2n}}{1+4 T_{2n}^2M_{2n}^2}=+\infty\, ,    
\end{equation}
then
system $(\ref{Wo.1})-(\ref{Wo.3})$ is asymptotically stable, that is any solution  $u$ of $(\ref{Wo.1})-(\ref{Wo.3})$ satisfies 
$E(u,t)\rightarrow 0$ for $t\rightarrow +\infty\,.$ 
\end{Theorem}

\begin{Theorem}\label{expo}
Assume   $\omega_1=\omega_2$, 
$\mbox{\rm i}_w), \mbox{\rm ii}_w)$ and $(\ref{new})$ are satisfied. 
Assume also
 that the observability inequality $(\ref{OCo})$ holds for  every time
$T>\overline T$
and that
$$
T_{2n}=T^*\quad\forall\ n\in\nat ,$$
with $T^*$ satisfying $(\ref{quartaA}),$
and
$$
T_{2n+1}=\tilde T\quad\forall\ n\in\nat .$$

Moreover,  assume that
\begin{equation}\label{ASS1Ao}
\sup_{n\in\nat} \ e^{(\xi +\frac 1 {\xi})M_{2n+1}\tilde T}c_n=d<1,
\end{equation}
where $c_n$ is as in $(\ref{no})$. Then, there exist two positive constants $\gamma,\mu$ 
such that
\begin{equation}\label{expestimateAo}
E(t)\le \gamma e^{-\mu t} E(0),\ \ t>0,
\end{equation}
for any solution of problem  $(\ref{Wo.1})-(\ref{Wo.3}).$
\end{Theorem}

\begin{Remark}{\rm
The case $\omega_1=\omega_2$ was already considered in \cite{ADE2012}.
Note that we significantly improve previous stability results (cfr. Theorems 
4.3 and 4.4 in \cite{ADE2012}).
}\end{Remark}

Under the restriction (\ref{rest}) we can obtain stability results 
also in the case $\omega_1\neq\omega_2.$
Note that the sets $\omega_1$ and $\omega_2$ may also have empty intersection. 

Indeed, we can restate Proposition \ref{intOK2abstraitnew}
and Theorem \ref{stab2abstraitnew} for our concrete model.

\begin{Proposition}\label{intOK2abstraitnewo}
Assume $\mbox{\rm i}_w), \mbox{\rm  ii}_w)$ and $(\ref{rest}).$
Moreover,
we assume that the observability inequality $(\ref{OCo})$ holds for every  time $T>\overline T$
and that, denoting  $T^*:=\inf_n \{T_{2n}\},$ 
$(\ref{quartaA})$ is satisfied.
Then, for any solution of system $(\ref{Wo.1})-(\ref{Wo.3})$ we have
\begin{equation}\label{stimabuona2Anewo}
E_S(t_{2n+1})\le \hat c_n E_S(t_{2n}),\quad \forall\ n\in\nat,
\end{equation}
where
\begin{equation}\label{nnewo}
\displaystyle{
\hat c_n= 
\frac{2c(1+4T_{2n}^2M_{2n}^2) }
{m_{2n}+ 2c(1+4T_{2n}^2M_{2n}^2)}
},
\end{equation}
  $c$ being the observability constant in $(\ref{OCo})$ corresponding to the time $T^*.$
\end{Proposition}

\begin{Theorem}\label{stab2abstraitnewo} 
Under the assumptions of Proposition $\ref{intOK2abstraitnewo},$
if  
$(\ref{star})$ holds,
then
system $(\ref{Wo.1})-(\ref{Wo.3})$ is asymptotically stable, that is any solution  $u$ of $(\ref{Wo.1})-(\ref{Wo.3})$ satisfies 
$E_S(u,t)\rightarrow 0$ for $t\rightarrow +\infty\,.$ 
\end{Theorem}

Under more restrictive assumption also an exponential stability estimate 
holds.

\begin{Theorem}\label{expnewo}
Assume $\mbox{\rm i}_w), \mbox{\rm ii}_w)$ and $(\ref{rest}).$
Assume also
 that the observability inequality $(\ref{OCo})$ holds for  every time
$T>\overline T$
and that
$$
T_{2n}=T^*\quad\forall\ n\in\nat ,
$$
with $T^*$ satisfying $(\ref{quartaA}),$
and
$$
T_{2n+1}=\tilde T\quad\forall\ n\in\nat ,
$$
with $\tilde T\le\tau.$

Moreover,  assume that
\begin{equation}\label{ASS1Anewnew}
\sup_{n\in\nat} \left [ 
e^{M_{2n+1}\tilde T}(\hat c_n+1)-1
\right ]
 =\hat d<1,
\end{equation}
where $\hat c_n$ is as in $(\ref{nnewo})$. Then, there exist two positive constants $\hat\gamma,\hat \mu$ 
such that
\begin{equation}\label{expestimateAnewo}
E_S(t)\le \hat\gamma e^{-\hat \mu t} E_S(0),\ \ t>0,
\end{equation}
for any solution of problem  $(\ref{Wo.1})-(\ref{Wo.3}).$
\end{Theorem}

\begin{Remark}{\rm
Note that the case $\omega_1\neq\omega_2$ is not covered from the abstract
setting of \cite{ADE2012}. 
}\end{Remark}

\bigskip

 {\em E-mail address,}
\\
\quad Serge Nicaise: \quad{\tt \bf
snicaise@univ-valenciennes.fr}
\\
\quad Cristina Pignotti: \quad{\tt \bf
pignotti@univaq.it}

\end{document}